\documentclass[11pt]{amsart}

\usepackage{hyperref}
\usepackage{amsmath,latexsym,amssymb,amsthm,array,amsfonts}

\usepackage{mathrsfs,dsfont,bbm}

\numberwithin{equation}{section}

\newtheorem{definition}{Definition}[section]
\newtheorem{theorem}{Theorem}[section]
\newtheorem{lemma}{Lemma}[section]
\newtheorem{corollary}{Corollary}[section]
\newtheorem{proposition}{Proposition}[section]
\newtheorem{example}{Example}[section]

\numberwithin{equation}{section}
\begin{document}

\title[Integral wrt the $G$-Brownian local time]{Integral
with respect to the $G$-Brownian local time}

\footnote[0]{${}^{*}$The Project-sponsored by NSFC (11171062), and Innovation Program of Shanghai Municipal Education Commission
(12ZZ063).

${}^{\S}$Corresponding author (litanyan@dhu.edu.cn).}

\author[L. Yan, X. Sun and B. Gao]
{Litan Yan${}^{\S}$, Xichao Sun and Bo Gao}

\date{}

\maketitle

\begin{center}
{\it Department of Mathematics, Donghua University, 2999 North
Renmin Rd., Songjiang, Shanghai 201620, P.R. China}
\end{center}

%%%%%%%%%%%%%%%%%%%%%
\maketitle

\begin{abstract}
Let ${\mathscr L}$ be the local time of $G$-Brownian motion $B$. In this paper, we prove the existence of the quadratic covariation $\langle f(B),B\rangle_{t}$ and the integral $\int_{\mathbb R}f(x){\mathscr L}(dx,t)$. Moreover, a sublinear version of the Bouleau-Yor identity
$$
\int_{\mathbb R}f(x){\mathscr L}(dx,t)=-\langle f(B),B\rangle_{t}
$$
is showed to hold under some suitable conditions. These allow us to write the It\^o's formula for $C^1$-functions.

\hfill
\begin{itemize}
\item  {\bf  Key words and phrases}: Nonlinear
expectation, $G$-Brownian motion, Local times, Young integration,
Quadratic covariation;

\item  {\bf 2010 Mathematics Subject Classification}: 60G05, 60G20, 60H05.
\end{itemize}
\end{abstract}

\section{Introduction}\label{sec1}
Motivated by various types of uncertainty and financial problems,
Peng~\cite{Peng1} has introduced a new notion of nonlinear
expectation, the so-called {\em $G$-expectation} (see also
Peng~\cite{Peng2,Peng3,Peng4}), which is associated with the
following nonlinear heat equation
$$
\begin{cases}
&\frac{\partial }{\partial t}u(t,x)=G(\Delta u),\qquad (t,x)
\in [0,+\infty)\times {\mathbb R},\\
&u(0,x)=\varphi(x),
\end{cases}
$$
where $\Delta$ is Laplacian, the sublinear function $G$ is defined
as
$$
G(\alpha)=\frac12\left(\overline{\sigma}^2\alpha^{+}
-\underline{\sigma}^2\alpha^{-}\right),\quad \alpha\in {\mathbb R}
$$
with two given constants $0<\underline{\sigma}<\overline{\sigma}$.
Together with the notion of $G$-expectations Peng also introduced
the related {\em $G$-normal distribution}, the {\em $G$-Brownian
motion} and related {\em stochastic calculus under $G$-expectation},
and moreover an It\^o's formula for the $G$-Brownian motion was
established. $G$-Brownian motion has a very rich and interesting new
structure which non-trivially generalizes the classical one. Briefly
speaking, a $G$-Brownian motion $B$ is a continuous process with
independent stationary increments $B_{t+s}-B_{t}$ being $G$-normally
distributed under a given {\em sublinear expectation} $\hat{{\mathbb
E}}$. A very interesting new phenomenon of $G$-Brownian motion $B$
is that its quadratic process $\langle B\rangle$ is a continuous
process with independent and stationary increments, but not a
deterministic process.

On the other hand, in the theory and applications of classical
stochastic calculus, the It\^o's formula plays a central role. But,
the restriction of It\^o's formula to functions with twice
differentiability often encounter difficulties in applications. For
the classical It\^o formula, many authors have given some extension
to the It\^o formula. Some approach extending It\^{o}'s formula are
to use the classical quadratic covariation and the {\em local
time-space calculus}. More works for the problem can be found in
Bouleau--Yor~\cite{Boul}, Eisenbaum~\cite{Eisen1,Eisen2}, Elworthy
{\em et al.}~\cite{Elworthy}, Feng-Zhao~\cite{Feng,Feng3}, F\"ollmer
{\it et al.}~\cite{Follmer}, Peskir~\cite{Peskir1},
Russo--Vallois~\cite{Russo2}, Yan--Yang~\cite{Yan4}, and the
references therein.

Recently, Li-Peng~\cite{Li-Peng} gave a more general It\^o integral
with respect to $G$-Brownian motion and It\^o's formula. Moreover,
Lin~\cite{Lin} studied local time of $G$-Brownian motion and
established a Tanaka formula under $G$-expectation. These motivate
us to consider the extensions of the classical {\em local time-space
calculus} under $G$-expectation. Let us first recall some known
results concerning the quadratic variation and It\^o's formula. Let
$F$ be an absolutely continuous function with locally square
integrable derivative $f$, that is,
$$
F(x)=F(0)+\int_0^xf(y)dy
$$
with $f$ being locally square integrable. Bouleau-Yor~\cite{Boul}
and F\"ollmer {\it et al}~\cite{Follmer} introduced the following
formulas:
\begin{equation}\label{sec1-eq1.1-1}
\left[f(W),W\right]_t=-\int_{\mathbb R}f(x){\mathscr L}^W(dx,t)
\end{equation}
and
\begin{equation}\label{sec1-eq1.1}
F(W_t)=F(0)+\int_0^tf(W_s)dW_s+\frac12\left[f(W),W\right]_t,
\end{equation}
where $W$ is the classical Brownian motion, ${\mathscr L}$ is local
time of Brownian motion $W$ and $\left[f(W),W\right]$ is the
classical quadratic covariation of $f(W)$ and $W$. If $f\in
C^1({\mathbb R})$,~\eqref{sec1-eq1.1} is the classical It\^o
formula. This result has been extended to some classical
semimartingales, smooth nondegenerate martingales and fractional
Brownian motion by Russo--Vallois~\cite{Russo2},
Moret--Nualart~\cite{Moret} and Yan {\it et al.}~\cite{Yan1,Yan3}.

In this paper we consider the identity~\eqref{sec1-eq1.1-1} and
It\^o's formula~\eqref{sec1-eq1.1} under {\em sublinear
expectation}. Our start point is to define the integral
\begin{equation}\label{sec1-eq1.2}
\int_{\mathbb R}f(x){\mathscr L}(dx,t)
\end{equation}
by using Young integration, where ${\mathscr L}$ is local time of
$G$-Brownian motion. If $f$ is of bounded $p$-variation with $1\leq
p<2$ we show that the integral~\eqref{sec1-eq1.2} and quadratic
covariation $\langle f(B), B\rangle$ exist, where the quadratic
covariation $\langle f(B),B\rangle$ is defined as
\begin{align*}
\langle f(B),B\rangle_{t}:=\lim_{n\rightarrow\infty}
\sum_{k=0}^{n-1}\{f(B_{t_{k+1}})-f(B_{t_{k}})\}(B_{t_{k+1}}
-B_{t_{k}})
\end{align*}
in ${\mathbb L}^{1}(\Omega)$, where $\{t_{k}\}$ is a partitions of
$[0,T]$ such that $\max\limits_{k}\{t_{k}-t_{k-1}\}\to 0$ as $n\to
\infty$.

This paper is organized as follows. In Section~\ref{sec2} we present
some preliminaries for $G$-Brownian motion. In Section~\ref{sec3},
we show that the one parameter integral~\eqref{sec1-eq1.2} exists
and establish the generalized It\^o formula
\begin{align}
F(B_{t})=F(0)+\int_{0}^{t}f(B_{s})dB_{s}-\frac{1}{2}\int_{\mathbb R}
f(x){\mathscr L}(dx,t),
\end{align}
where $F\in C^1({\mathbb R})$ and $F'=f$ is of bounded $p$-variation
with $1\leq p<2$. In Section~\ref{sec4} we prove the existence of
the quadratic covariation $\langle f(B), B\rangle$. As a result, we
get a sublinear version of the Bouleau-Yor identity
\begin{equation}\label{sec1-eq1.3}
\langle f(B),B\rangle_t=-\int_{\mathbb {R}}f(x){\mathscr L}(dx,t),
\end{equation}
under some suitable conditions. In Section~\ref{sec5} we extend
these results to the time-dependent case. In Appendix we give the
other representation of quadratic covariation:
\begin{align*}
\frac1{\varepsilon}\int_0^t
(B_{s+\varepsilon}-B_{s})^2ds\longrightarrow \langle B\rangle_{t}
\end{align*}
in ${\mathbb L}^1$, as $\varepsilon\downarrow 0$.

%%%%%%%%%%%%%%%%%%%%%%%%%%%%%%%%%%%%%%%%%%%%%%%%%%%%%%%%%%%%%%%%%%%%%
%%%%%%%%%%%%%%%%%%%%%%%%%%%%%%%%%%%%%%%%%%%%%%%%%%%%%%%%%%%%%%%%%%%%%
%%%%%%%%%%%%%%%%%%%%%%%%%%%%%%%%%%%%%%%%%%%%%%%%%%%%%%%%%%%%%%%%%%%%%

\section{Preliminaries}\label{sec2}
In this section, we briefly recall some basic notations and results
for $G$-Brownian motion under $G$-framework. For more aspects on
these material we refer to Li-Peng~\cite{Li-Peng}, Lin~\cite{Lin}
and Peng~\cite{Peng1,Peng2,Peng3,Peng4}. For simplicity throughout
this paper we let $C$ stand for a positive constant depending only
on the subscripts and its value may be different in different
appearance.

\subsection{Sublinear expectation space}

Let $\Omega\neq \emptyset$ be a given set and let $\mathcal {H}$ be
a linear space of real valued functions defined on  $\Omega$ such
that $1\in \mathcal {H}$ and $|X|\in \mathcal {H}$ for all $X\in
\mathcal {H}$.

\begin{definition}\label{def2-1}
A sublinear expectation $\hat{\mathbb{E}}$ on $\mathcal {H}$ is a
functional with the following properties : for all $X,Y\in \mathcal
{H}$, we have
\begin{itemize}
\item {\sc {Monotonicity} : } if $X\geq Y$, then $\hat{\mathbb{E}}[X]
\geq\hat{\mathbb{E}}[Y]$;
\item {\sc {Constant preserving} :} $\hat{\mathbb{E}}[c]=c $,
for all~$c\in \mathbb{R}$;
\item {\sc {Sub-additivity} :}$\hat{\mathbb{E}}[X]- \hat{\mathbb{E}}
[Y]\leq\hat{\mathbb{E}}[X- Y]$;
\item {\sc {Positive homogeneity} :} $\hat{\mathbb{E}}[\lambda X]
=\lambda \hat{\mathbb{E}}[X]$, for all $\lambda\geq 0$.
\end{itemize}
The triple $(\Omega,{\mathcal H},\hat{\mathbb{E}})$ is called a
sublinear expectation space, and ${\mathcal H}$ is considered as the
space of random variables on $\Omega$.
\end{definition}
It is important to note that we can suppose that
$$
\varphi(X_1,\ldots,X_d)\in {\mathcal H}
$$
if $X_i\in {\mathcal H}, i = 1,\ldots, d$, for all $\varphi\in
C_{b,Lip} ({\mathbb R}^d)$, where $C_{b,Lip}({\mathbb R}^d)$ denotes
the space of all bounded and Lipschitz functions on ${\mathbb R}^d$.
In a sublinear expectation space $(\Omega,{\mathcal H},
\hat{\mathbb{E}})$, a random vector $Y=(Y_1,\ldots,Y_n),Y_i\in
{\mathcal H}$ is said to be independent under $\hat{\mathbb{E}}$
from another random vector $X=(X_1,\ldots,X_m),X_i\in {\mathcal H}$,
if for each test function $\varphi\in C_{b,Lip}({\mathbb R}^{m+n})$
we have
$$
\hat{\mathbb{E}}[\varphi(X,Y)]=\hat{\mathbb{E}}\left[
\hat{\mathbb{E}}[\varphi(x,Y)]_{x=X}\right].
$$
Two $n$-dimensional random vectors $X$ and $Y$ defined respectively
in the sublinear expectation spaces $(\Omega_1,{\mathcal H}_1,
\hat{\mathbb{E}}_1)$ and $(\Omega_2,{\mathcal H}_2,
\hat{\mathbb{E}}_2)$ are called identically distributed, denoted by
$X\sim Y$, if
$$
\hat{\mathbb{E}}_1[\varphi(X)]=\hat{\mathbb{E}}_2[\varphi(Y)]
$$
for all $\varphi\in C_{b,Lip}({\mathbb R}^n)$.

Let $\underline{\sigma},\overline{\sigma}$ be two real numbers with
$0<\underline{\sigma}<\overline{\sigma}$. A random variable $X$ in a
sublinear expectation space $(\Omega,{\mathcal H},\hat{\mathbb{E}})$
is called $G$-normal distributed, denoted by $\xi\sim N(0,
[\underline{\sigma}^2,\overline{\sigma}^2])$, if for each
$\varphi\in C_{b,lip}({\mathbb R})$, the function defined by
$$
u(t, x) :=\hat{\mathbb{E}}\left[\varphi(x+\sqrt{t}\xi\right],\quad
(t,x)\in [0,\infty)\times {\mathbb R}
$$
is the unique viscosity solution of the following nonlinear heat
equation :
$$
\begin{cases}
&\frac{\partial }{\partial t}u(t,x)=G(\Delta u),\qquad (t,x)\in
[0,+\infty)\times {\mathbb R},\\
&u(0,x)=\varphi(x),
\end{cases}
$$
where $\Delta$ is Laplacian and the sublinear function $G$ is
defined as
$$
G(\alpha)=\frac12\left(\overline{\sigma}^2\alpha^{+} -\underline{
\sigma}^2\alpha^{-}\right),\quad \alpha\in {\mathbb R}.
$$
\begin{example}[Peng~\cite{Peng1}]
Let $\xi\sim N(0,[\underline{\sigma}^2,\overline{\sigma}^2])$. We
then have
$$
\hat{{\mathbb E}}\left[\varphi(\xi)\right]=\frac{1}{\sqrt{2\pi}
\overline{\sigma}} \int_{\mathbb
R}\varphi(x)e^{-\frac1{2\overline{\sigma}^2}x^2}dx
$$
for all convex functions $\varphi$, and
$$
\hat{{\mathbb E}}\left[\psi(\xi)\right]=\frac{1}{\sqrt{2\pi}
\underline{\sigma}} \int_{\mathbb
R}\varphi(x)e^{-\frac1{2\underline{\sigma}^2}x^2}dx
$$
for all concave functions $\psi$.
\end{example}
\subsection{Some spaces}
In this paper we throughout let $\Omega=C_{0}(\mathbb{R}^{+})$ be
the space of all real valued continuous functions on $[0,\infty)$
with initial value $0$, equipped with the distance
$$
\rho(\omega^1,\omega^2)=\sum_{i=1}^\infty 2^{-i}\left[\left(\max_{
t\in [0,i]}|\omega_t^1-\omega_t^2|\right)\wedge 1\right],\quad
\omega^1,\omega^2\in \Omega.
$$
We denote by ${\mathscr B}(\Omega)$ the Borel-algebra on $\Omega$.
We also denote, for each $t\in [0,\infty)$,
$$
\Omega_t=\{\omega_{\cdot\wedge t},\omega\in \Omega\},
$$
and ${\mathscr F}_t={\mathscr B}(\Omega_t)$, where $x\wedge y=
\min\{x,y\}$. We also denote

\begin{itemize}
\item $L^0(\Omega)$ : the space of all ${\mathscr B}(\Omega)$-measurable
real valued functions on $\Omega$;
\item $L^0(\Omega_t)$ : the space of all ${\mathscr B}(\Omega_t)$-measurable
real valued functions on $\Omega_t$;
\item $L_b(\Omega)$ : the space of all bounded elements in $L^0(\Omega)$;
\item $L_b(\Omega_t)$ : the space of all bounded elements in $L^0(\Omega_t)$.
\end{itemize}
Let ${\mathbb L}^p_G(\Omega)$ be the closure of ${\mathcal H}$ with
respect to the norm
$$
\|X\|_p=\hat{\mathbb{E}}\left[|X|^p\right]^{1/p}
$$
with $p\in [1,\infty)$. Clearly, the space ${\mathbb L}^p_G(\Omega)$
is a Banach space and the space $C_b(\Omega)$ of bounded continuous
functions on $\Omega$ is a subset of ${\mathbb L}^1_G(\Omega)$, and
moreover, for the sublinear expectation space $(\Omega,{\mathbb
L}^p_G(\Omega),\hat{\mathbb{E}})$  there exists a weakly compact
family ${\mathcal P}$ of probability measures on $(\Omega,{\mathscr
B}(\Omega))$ such that
$$
\hat{\mathbb{E}}=\sup_{P\in {\mathcal P}}E_P.
$$
So we can introduce the Choquet capacity $\hat{C}$ by taking
$$
\hat{C}(A)=\sup_{P\in {\mathcal P}}P(A),\qquad A\in{\mathscr B}(
\Omega).
$$
\begin{definition}
A set $ A\subset\Omega$ is called polar if $\hat{C}(A)=0$. A
property is said to hold "quasi surely" ({\rm q.s.}) if it holds
outside a polar set.
\end{definition}

By using the above family of probability measures $P$ we can
characterize the space ${\mathbb L}_G^p(\Omega)$ as
\begin{align*}
{\mathbb L}_G^p(\Omega)&=\Bigl\{L^0(\Omega)\ni X {\text { is
continuous, q.s., and }} \sup_{P\in {\mathcal
P}}E_P[|X|^p] <\infty\Bigr\}\\
&\equiv\Bigl\{L^0(\Omega)\ni X {\text { is continuous, q.s., and }}
\lim_{n\to \infty}\sup_{P\in {\mathcal P}}E_P[|X|^p1_{\{|X|>n\}}]
=0\Bigr\}.
\end{align*}
The following three results can be consulted in Denis {\em et
al}~\cite{Denis} and Hu-Peng~\cite{Hu-Peng}.
\begin{lemma}[Denis {\em et al}~\cite{Denis} and
Hu-Peng~\cite{Hu-Peng}]\label{lemma2.1000}
Let
$\{X_n,n=1,2,\ldots\}$ be an monotonically decreasing sequence of
nonnegative random variances in $C_b(\Omega)$. If $X_n$ converges to
zero {\em q.s} on $\Omega$, then we have
$$
\lim_{n\to 0}\hat{\mathbb{E}}[X_n]=0.
$$
Moreover, if $X_n\uparrow X$ and
$\hat{\mathbb{E}}[X],\hat{\mathbb{E}}[X_n]$ are finite for all
$n=1,2,\ldots$, we then have
$$
\lim_{n\to 0}\hat{\mathbb{E}}[X_n]=\hat{\mathbb{E}}[X].
$$
\end{lemma}
\begin{lemma}[Denis {\em et al}~\cite{Denis} and Hu-Peng~\cite{Hu-Peng}]
Let $1\leq p<\infty$. Consider the sets ${\mathbb L}_G^p(\Omega)$
and ${\mathbb L}^p={\mathcal L}^p/{\mathcal N}$, where
\begin{align*}
{\mathcal L}^p&=\Bigl\{X\in L^0(\Omega) :\;
\hat{\mathbb{E}}(|X|^p)=\sup_{P\in {\mathcal P}}E_P[|X|^p]<\infty\Bigr\},\\
{\mathcal N}&=\Bigl\{X\in L^0(\Omega) :\; X=0\;\;\;q.s.\Bigr\}.
\end{align*}
Then
\begin{itemize}
\item ${\mathbb L}^p$ is a Banach space with respect to the norm
$\|\cdot\|_p$;
\item ${\mathbb L}^p_G$ is the completion of $C_b(\Omega)$ with
respect to the norm $\|\cdot\|_p$.
\end{itemize}
\end{lemma}
\begin{lemma}[Denis {\em et al}~\cite{Denis} and Hu-Peng~\cite{Hu-Peng}]
\label{lem2.3} For a given $p\in (0, +\infty]$, If the sequence
${\mathbb L}^p\supset \{X_n\}$ converges to $X$ in ${\mathbb L}^p$,
then there exists a subsequence $\{X_{n_k}\}$ such that $X_{n_k}$
converges to $X$ quasi-surely.
\end{lemma}

We denote by ${\mathbb L}_\ast^p(\Omega)$ the completion of
$L_b(\Omega)$ with respect to the norm $\|\cdot\|_p$.

\subsection{$G$-Brownian motion}
Now, let us recall the definition of $G$-Brownian motion and related
It\^o's integral.
\begin{definition}
A process $B=\{B_t, t\geq 0\}\subset {\mathcal H}$ in a sublinear
expectation space $(\Omega,{\mathcal H},\hat{\mathbb{E}})$ is called
a $G$-Brownian motion if the following properties are satisfied:
\begin{itemize}
\item $B_0=0$;
\item For each $t,s\geq 0$, the increment $B_{t+s}-B_t$ is $N(0,
[\underline{\sigma}^2s,\overline{\sigma}^2s])$-distributed and is
independent from $(B_{t_1},\ldots,B_{t_n})$, for all
$n=0,1,2,\ldots$ and $0\leq t_1\leq t_2\leq \cdots\leq t_n\leq t$.
\end{itemize}
\end{definition}
The $G$-Brownian motion $B$ has following properties :
\begin{itemize}
\item [(1)] For all $\xi\in {\mathbb L}^2(\Omega_t)$, we have
$\hat{\mathbb{E}}[\xi(B_T-B_t)]=0$ with $0\leq t\leq T$;
\item [(2)] For all ${\mathscr B}(\Omega_t)$-measurable real valued,
bounded functions $\xi$, we have
$$
\hat{\mathbb{E}}[\xi^2(B_T-B_t)^2]\leq \overline{\sigma}^2(T-t)
\hat{\mathbb{E}}[\xi^2],\quad 0\leq t\leq T;
$$
\item [(3)] For all $t\geq 0$, we have $\hat{\mathbb{E}}[B_t]
=\hat{\mathbb{E}}[-B_t]=0$;
\item [(4)] $t\mapsto B_t$ is H{\"o}lder continuous of order
$\delta<\frac12$, quasi-surely.
\end{itemize}
In Li-Peng~\cite{Li-Peng}, a generalized It\^o integral and a
generalized It\^o formula with respect to the G-Brownian motion are
introduced. For arbitrarily fixed $p\geq 1$ and $T\in {\mathbb
R}_{+}$, we denote by $M_b^{p,0}([0,T])$ the set of step processes:
\begin{equation}\label{sec2.3-eq1}
\eta_{t}(\omega)=\sum_{j=1}^N\xi_j(\omega)1_{[t_{j-1},
t_{j})}(t),\qquad \xi_j \in L_b(\Omega_{t_j})
\end{equation}
with $0=t_{0}<\cdots<t_{N}=T$. For the process of the
form~\eqref{sec2.3-eq1} we define the related Bochner integral as
follows
$$
\int_0^T\eta_tdt=\sum_{j=1}^N\xi_j(\omega)(t_{j-1}, t_{j}).
$$
For every $\eta\in M_b^{p,0}([0,T])$ we set
$$
\hat{\mathbb E}_T(\eta):=\frac1T\hat{\mathbb E}\int_0^T\eta_tdt.
$$
Then $\hat{\mathbb E}_T$ forms a sublinear expectation. Moreover, we
denote by $M_\ast^{p}([0,T])$ the completion of $ M_b^{p,0}([0,T])$
under the norm
\begin{align*}
\|\eta\|_{M_\ast^{p}([0,T])}&=\left( \hat{\mathbb{E}}[\int_{0}^{T}|
\eta_{s}|^{p}ds] \right)^{\frac{1}{p}}.
\end{align*}

\begin{definition}
For every $\eta \in M_b^{p,0}([0,T])$ of the
form~\eqref{sec2.3-eq1}, we define the It\^o integral of $\eta$ with
respect to $G$-Brownian motion $B$ as
\begin{align*}
I(\eta):=\int_{0}^{T}\eta_sdB_s=\sum_{j=1}^N\xi_{j}
(B_{s_j}-B_{s_{j-1}}).
\end{align*}
\end{definition}
The mapping $I : M_b^{p,0}([0, T])\to {\mathbb L}^2_\ast(\Omega_T)$
is a linear continuous mapping and thus can be continuously extended
to $I : M_\ast^2([0, T])\to {\mathbb L}^2_\ast(\Omega_T)$, which is
called the It\^o integral of $\eta\in M_\ast^2([0, T])$ with respect
to $G$-Brownian motion $B$, and define
$$
\int_{0}^{t}\eta_sdB_s=\int_0^T1_{\{0\leq s\leq t\}}\eta_sdB_s
$$
for all $\eta\in M_\ast^{2}([0,T])$ and $t\in [0,T]$. We have
$$
\hat{\mathbb{E}}\left(\int_0^T\eta_sdB_s\right)=0
$$
and
$$
\hat{\mathbb{E}}\left[\left(\int_0^T\eta_sdB_s\right)^2\right] \leq
\overline{\sigma}^2 \hat{\mathbb{E}}\left[\int_0^T\eta_s^2ds\right]
$$
for all $\eta\in M_\ast^{2}([0,T])$. Moreover, the process
$\{\int_{0}^{t}\eta_sdB_s,t\in [0,T]\}$ is continuous in $t$ quasi
surely, and
$$
\int_{0}^{\cdot}\eta_sdB_s\in M_\ast^{2}([0,T])
$$
for all $\eta\in M_\ast^{2}([0,T])$.

\begin{definition}[Quadratic Variation]
Let $\pi _{t}^{N}=\{0=t_{0}^{N}<t_{1}^{N}<\cdots<t_{N-1}^{N}=t\}$ be
a partition of $[0,t]$ for $t>0$, such that $\mu(\pi
_{t}^{N}):=\max\limits_{j}\{t_j-t_{j-1}\}\to 0$ as
$N\rightarrow\infty$. The quadratic variation of $G$-Brownian motion
$B$ is defined as
\begin{align*}
\langle B\rangle_{t}=\lim_{\mu(\pi _{t}^{N})\to
0}\sum_{k=0}^{N-1}(B_{t_{k+1}^{N}}
-B_{t_{k}^{N}})^{2}=B_t^2-2\int_0^tB_sdB_s
\end{align*}
in ${\mathbb L}_{G}^{2}(\Omega)$.
\end{definition}
The function $t\mapsto \langle B\rangle_t$ is continuous and
increasing outside a polar set. We can define the integral
\begin{align*}
\int_{0}^{T}\eta_{t}d\langle B\rangle_{t}:=\sum_{j=1}^N\xi_j(\langle
B\rangle_{t_j}-\langle B\rangle_{t_{j-1}})
\end{align*}
as a map from $M_b^{1,0}([0,T])$ into ${\mathbb
L}_\ast^{1}(\Omega_T)$, and the map is linear and continuous, and it
can be extended continuously to $M_\ast^{1}([0,T])$.

\begin{theorem}[It\^o's formula]
Let $F\in C^{2}(\mathbb{R}\times {\mathbb R}_+)$. We then have
\begin{align*}
F(B_{t},t)=F(B_{0},0)&+\int_{0}^{t}\frac{\partial}{\partial
x}F(B_{s},s)dB_{s}\\
&+\int_{0}^{t}\frac{\partial}{\partial
t}F(B_{s},s)ds+\frac{1}{2}\int_{0}^{t}\frac{\partial^2 }{\partial
x^{2}}F(B_{s},s)d\langle B \rangle_{s}
\end{align*}
for all $t\geq 0$.
\end{theorem}
Finally, recall that the $G$-Brownian motion $B$ has a jointly
continuous local time ${\mathscr L}(x,t)$ which satisfies (see
Lin~\cite{Lin})
\begin{align*}
{\mathscr L}(x,t)=\lim_{\varepsilon\rightarrow
0}\frac{1}{2\varepsilon}\int_{0}^{t}1_{(x-\varepsilon,x+\varepsilon)}
(B_{s})d\langle B\rangle_{s}
\end{align*}
in ${\mathbb L}^{2}$, and $(x,t)\mapsto {\mathscr L}(x,t)$ is
H\^older continuous of order $0<\gamma<\frac12$. Moreover, the
following Tanaka formula holds:
\begin{align*}
|B_{t}-x|=|x|+\int_{0}^{t}{\rm sign}(B_{s}-x)dB_{s}+{\mathscr
L}(x,t).
\end{align*}
for all $x\in\mathbb{R}$ and $t\geq 0$. As a result of Tanaka formula we have
$$
{\mathscr L}(x,t)=0
$$
if $|x|\geq \sup\limits_{0\leq s\leq t}|B_s|$, i.e., the function
$x\mapsto {\mathscr L}(x,t)$ has a compact support for all $0\leq
t<\infty$.
\begin{theorem}[Theorem 5.4 in Lin~\cite{Lin}]
Let $\underline{\sigma}> 0$ and let $p\geq 1,a\leq b$. Then the
convergence in ${\mathbb L}^p$
$$
\lim_{n\to \infty}\sum_{i=0}^{2^n-1}\left({\mathscr
L}(a^n_{i+1},t)-{\mathscr L}(a^n_i,t)\right)^2=4\int_a^b{\mathscr
L}(x,t)dx
$$
holds uniformly in $t\in [0, T ]$ for all the sequence of partitions
$\pi_n=\{a_i^n = a + 2^{-n}i(b-a),i=0, 1, 2,\ldots,2^n\},\;n\geq 1$,
of the interval $[a, b]$.
\end{theorem}

%%%%%%%%%%%%%%%%%%%%%%%%%%%%%%%%%%%%%%%%%%%%%%%%%%%%%%%%%%%%%%%%%%
%%%%%%%%%%%%%%%%%%%%%%%%%%%%%%%%%%%%%%%%%%%%%%%%%%%%%%%%%%%%%%%%%%

\section{One Parameter Integrals of Local Time}\label{sec3}
In this section, we study the integral
\begin{equation}\label{sec3-eq3.1}
\int_{\mathbb{R}}f(x){\mathscr L}(dx,t)
\end{equation}
where $f$ is a real function, $B$ is a $G$-Brownian motion and
${\mathscr L}$ is the local time of $B$. To use Young integration to
establish the integral, we first investigate $p$-variation of the
mapping $x\mapsto {\mathscr L}(x,t)$ for every $t\geq 0$. Throughout
this paper we let $\underline{\sigma}>0$.

\begin{lemma}[Occupation times formula]\label{lem3.1}
For all $t\geq 0$ and every bounded function $f$, we have
\begin{align*}
\int_{0}^{t}f(B_{s})d\langle B \rangle_{s}=\int_{\mathbb
R}f(x){\mathscr L}(x,t)dx.
\end{align*}
\end{lemma}

This formula is first introduced by Lin~\cite{Lin} for
$f(x)=1_{[a,b)}(x)$ (see Theorem 5.2 in Lin~\cite{Lin}). By
approximation we can give the above formula. In fact, for any simple
function $f_\Delta(x)=\sum_{i} a_{i}1_{[x_{i-1},x_{i})}(x)$, where
$\{x_{0}, x_{1},\cdots\}$ is a partition of ${\mathbb R}$, we have
\begin{align*}
\int_{0}^{t}f_{\Delta}(B_{s})d{\langle B \rangle_{s}}=\int_{\mathbb
R}f_\Delta(x){\mathscr L}(x,t)dx ,\quad q.s.
\end{align*}
by the linearity of the integral.

Now, for every bounded function $f$, one can show that there is a
sequence of bounded simple functions $f_n,\;n\geq 1$ such that
$f_n(x)\uparrow f(x)$ for all $x\in {\mathbb R}$ (of course, there
also exists a sequence of simple functions conversing to $f$
uniformly)
$$
\int_{0}^{t}f_n(B_{s})d{\langle B \rangle_{s}}\rightarrow
\int_{0}^{t}f(B_{s})d{\langle B \rangle_{s}}
$$
and
$$
\int_{\mathbb R}f_{n}(x){\mathscr L}(x,t)dx\rightarrow \int_{\mathbb
R}f(x){\mathscr L}(x,t)dx,
$$
in ${\mathbb L}^1$, by Lemma~\ref{lemma2.1000}, as $n\to \infty$,
which gives
\begin{align*}
\int_{0}^{t}f(B_{s})d{\langle B \rangle_{s}}=\int_{\mathbb
R}f(x){\mathscr L}(x,t)dx,\quad q.s.
\end{align*}
and the lemma follows.

\begin{definition}
Let $p\geq 1$ be a fixed real number. A function $f:
[a,b]\mapsto\mathbb{R}$ is said to be of bounded p-variation if
\begin{align*}
v_p(f)=\sup_{\Delta_{n}}\sum_{i=0}^{n}|f(x_{i+1})-f(x_{i})|^{p} <\infty,
\end{align*}
where the supremum is taken over all partition
$\Delta_{n}=\{a=x_{0}<x_{1}<\cdot\cdot\cdot<x_{n}=b\}$ of $[a,b]$.
\end{definition}

\begin{corollary}\label{cor3.2}
Let $p>2$. Then the mapping $x\mapsto {\mathscr L}(x,t)$ is of bounded $p$-variation q.s. for any $0\leq t\leq T$.
\end{corollary}

We now can establish one parameter integral~\eqref{sec3-eq3.1}. Denote by ${\mathcal W}_p([a,b])$ ($p\geq 1$) the
set of all measurable functions $f$ on $[a,b]$ such that $
v_p(f)<\infty$. For $1\leq p<\infty$ define
$$
\|f\|_{(p)}:=v_p(f)^{1/p}.
$$
Then $\|\cdot\|_{(p)}$ is a seminorm on ${\mathcal W}_p([a,b])$,
which is called the $p$-variation seminorm. For $1\leq p\leq\infty$
define
$$
\|f\|_{[p]}:=\|f\|_{(p)}+\|f\|_{\infty},
$$
where $\|f\|_{\infty}=\sup_{x\in [a,b]}|f(x)|$. Then
$\|\cdot\|_{[p]}$ is a norm on ${\mathcal W}_p([a,b])$, which is
called the $p$-variation norm. The space $({\mathcal
W}_p([a,b]),\|\cdot\|_{[p]})$ is a Banach space for $p\geq 1$, and
$f\in {\mathcal W}_p$ means that $f\in {\mathcal W}_p([a,b])$ for
any $a,b\in {\mathbb R}$, and moreover $f\in {\mathcal W}_p$ is
locally bounded. For these, see Dudley-Norvai\v{s}a~\cite{Dudley}
and Young~\cite{Young}.
\begin{lemma}[Dudley-Norvai\v{s}a~\cite{Dudley}]
Let $f\in {\mathcal W}_p([a,b])$ and $g\in {\mathcal W}_q([a,b])$,
where $p,q\geq 1$ and $\frac1p+\frac1q>1$. If $f$ and $g$ have no
common discontinuities, then the Young integral (Bochner integral)
$$
\int_a^bf(x)dg(x):=\lim_{|\triangle_n|\rightarrow 0}\sum_{j=1}^{n}
f(\xi_{j})(g(x_{j})-g(x_{j-1}))
$$
exists, where $\xi_j\in[x_{j-1}, x_j]$ ($j=1,2,\ldots$),
$|\triangle_n|=\max_{1\leq j\leq n}|x_j-x_{j-1}|$ and the Love-Young
inequality
$$
\left|\int_a^bf(x)dg(x)\right|\leq C_{p,q}\|f\|_{[p]}\|g\|_{(p)}
$$
holds.
\end{lemma}

\begin{definition}
Let $f$ be a Borel function. We define Bochner integral
$$
\int_{\mathbb R}f(x){\mathscr L}(dx,t)=\lim_{|\Delta_{n}|\to 0}
\sum_{j=1}^{n}f(x_{j})\big({\mathscr L}(x_{j},t)-{\mathscr
L}(x_{j-1},t)\big)
$$
exists in ${\mathbb L}^1$ for any $0\leq t\leq T$, where
$\{a=x_{0}<x_{1}<\cdot\cdot\cdot<x_{n}=b\}$ is any partition of
$[a,b]$ and $|\Delta_{n}|=\max\{x_i-x_{i-1}\}$.
\end{definition}
\begin{proposition}\label{prop3.1}
Let $f$ be of bounded $p$-variation with $1\leq p<2$. Then
$$
\int_{\mathbb R}{\mathscr L}(x,t)df(x)
$$
exists in ${\mathbb L}^1$ and we have
$$
\int_{\mathbb R}f(x){\mathscr L}(dx,t)=-\int_{\mathbb R}{\mathscr
L}(x,t)df(x)
$$
for all $t\geq 0$. In particular, for $f\in C^{1}({\mathbb R})$ we
have
\begin{align*}
\int_{\mathbb R}f(x){\mathscr L}(dx,t)=-\int_{\mathbb R}{\mathscr
L}(x,t)f'(x)dx.
\end{align*}
\end{proposition}
\begin{proof}
Clearly, the Bochner integral
$$
\int_{\mathbb R}{\mathscr L}(x,t)df(x)
$$
exists q.s. for all $t\geq 0$, and
$$
\int_{\mathbb R}{\mathscr L}(x,t)df(x)\in {\mathbb L}^1,
$$
because $x\mapsto {\mathscr L}(x,t)$ is continuous and has a compact
support for each $0\leq t<\infty$. It follows that
\begin{align*}
\sum_{j=1}^{n}f(x_{j-1})\big({\mathscr L}(x_{j},t)
&-{\mathscr L}(x_{j-1},t)\big)\\
&=\sum_{j=1}^{n}f(x_{j-1}){\mathscr
L}(x_{j},t)-\sum_{j=0}^{n-1}f(x_{j}) {\mathscr
L}(x_{j},t)\\
&=-\sum_{j=1}^{n}(f(x_{j})-f(x_{j-1})){\mathscr L}(x_{j},t)
\end{align*}
by adding some points in the partition $\{x_i\}$ to make ${\mathscr
L}(x_{n},t)=0$ and ${\mathscr L}(x_{1},{t})=0$. This completes the
proof.
\end{proof}
Define the mollifier $\theta$ by
\begin{equation}\label{sec3-eq3.7}
\theta(x)=
\begin{cases}
c\exp(\frac{1}{(x-1)^2-1}), & {\text {$x\in(0,2)$}},\\
0, & {\text {$x\not\in(0,2)$}},
\end{cases}
\end{equation}
where $c$ is a normalizing constant such that $\int_{\mathbb
R}\theta(x)dx=1$. Set $\theta_n(x)=n\theta(nx)$. For a locally
integrable function $g(x)$ we define
$$
g_n(x)=\int_{\mathbb
R}\theta_n(x-y)g(y)dy=\int_0^2\theta(z)g(x-\frac{z}n)dz,\quad n\geq
1.
$$

\begin{lemma}\label{lem3.3}
Let $g$ be of bounded $p$-variation with $1\leq p<2$. Suppose that
$g_n$ is defined as above, then $g_n$ is of bounded $p$-variation
with $1\leq p<2$ and $g_n\in C^\infty({\mathbb R})$ for every $n\geq
1$, and moreover, the convergence
$$
\lim_{n\rightarrow \infty}\int_{\mathbb R}g_n(x){\mathscr L}(dx,t)=\int_{\mathbb R}g(x){\mathscr L}(dx,t)
$$
holds in ${\mathbb L}^1$.
\end{lemma}
\begin{proof}
We have
\begin{align*}
\left|\int_{\mathbb R}{\mathscr L} (x,t)dg_n(x)-\int_{\mathbb
R}{\mathscr L}(x,t)dg(x)\right|&=\left|\int_{\mathbb R}{\mathscr
L}(x,t)d(g_n(x)-g(x))\right|\\
&\leq C_{p,q}\|{\mathscr L}(\cdot,t)\|_{[p]}\|g_n(x)-g(x)\|_{(p)}
\end{align*}
for all $0\leq t<\infty$ because $x\mapsto {\mathscr L}(x,t)$ is
continuous and has a compact support for each $0\leq t<\infty$, and
the lemma follows.
\end{proof}
By using the above lemmas, we immediately get an extension of
It\^{o} formula stated as follows, which is an analogue of the Bouleau-Yor type formula.
\begin{theorem}\label{th3.1}
Let the function $f$ be of bounded $p$-variation with $1\leq p<2$
and let $F$ be an absolutely continuous function with derivative
$F'=f$. Then
\begin{equation}\label{th3.1-eq1}
F(B_{t})=F(0)+\int_{0}^{t}f(B_{s})dB_{s}-\frac{1}{2}\int_{\mathbb
R}f(x){\mathscr L}(dx,t).
\end{equation}
\end{theorem}

Recall that if $F$ is the difference of two convex functions, then $F$ is an absolutely continuous function with derivative of bounded variation. Thus, the It\^o-Tanaka formula
\begin{align*}
F(B_t)&=F(0)+\int_0^tF^{'}(B_s)dB_s+\frac12\int_{\mathbb
R}{\mathscr L}(x,t)F''(dx)\\
&\equiv F(0)+\int_0^tF^{'}(B_s)dB_s-\frac12\int_{\mathbb
R}F'(x){\mathscr L}(dx,t)
\end{align*}
holds.

\begin{proof}[Proof of Theorem~\ref{th3.1}]
For $n\geq1$, we set
\begin{align*}
F_{n}(x):=\int_{\mathbb R}\theta_{n}(y)F(x-y)dy\equiv
\int_0^2F(x-\frac{y}n)\theta(y)dy,
\end{align*}
where $\theta_{n},n\geq 1$ are the mollifiers defined as above. Then
$F_{n}\in C^{\infty}(\mathbb{R})$,
\begin{align*}
F_{n}(B_{t})=F_{n}(B_{0})+\int_{0}^{t}F'_{n}(B_{s})dB_{s}
+\frac{1}{2}\int_{0}^{t}F^{''}_{n}(B_{s})d{\langle B \rangle _{s}}
\end{align*}
and $F'_{n}$ is of bounded $p$-variation with $1\leq p<2$ for all
$n=1,2,\ldots$. Now, we prove the theorem in two steps.

{\bf Step I.} Let $f$ be bounded uniformly in ${\mathbb R}$. Then
$F_n\to F$ uniformly, as $n$ tends to infinity, because
\begin{align*}
|F_{n}(x)-F(x)|&\leq \int_0^2|F(x-\frac{y}n)-F(x)|\theta(y)dy\\
&\leq \sup_{z}|f(z)|\left(\int_0^2y\theta(y)dy\right)\frac1n
\end{align*}
for a constant $C>0$. Moreover, we have $F'_n,f\in {\mathcal W}_p$
for all $n=1,2,\ldots$, and $F'_n\to f$ in ${\mathcal W}_p$, as $n$
tends to infinity. Noting that
$$
|F'_n(x)-f(x)|\leq \|F'_n-f\|_{[p]},\qquad n=1,2,\ldots
$$
for all $x\in {\mathbb R}$, we get $F'_n\to f$ uniformly, as $n$
tends to infinity. It follows that
\begin{align*}
F_{n}(B_{t})\longrightarrow F(B_{t}),\qquad
\int_{0}^{t}F'_{n}(B_{s})dB_{s}\longrightarrow \int_{0}^{t}
f(B_{s})dB_{s},
\end{align*}
in ${\mathbb L}^1$, as $n\to \infty$. On the other hand, we have
\begin{align*}
\int_{0}^{t}F^{''}_{n}(B_{s})d{\langle B \rangle
_{s}}&=\int_{\mathbb R}{\mathscr L}(x,t)F^{''}_{n}(x)dx
=-\int_{\mathbb
R}F'_{n}(x){\mathscr L}(dx,t)\\
&\longrightarrow-\int_{\mathbb R}f(x){\mathscr L}(dx,t)
\end{align*}
in ${\mathbb L}^1$, by Lemma~\ref{lem3.3} and
Proposition~\ref{prop3.1}, which deduces, as $n\to \infty$
\begin{align*}
\int_{0}^{t}F'_{n}(B_{s})dB_{s}&=F_{n}(B_{t})-F_{n}(0)-
\frac{1}{2}\int_{0}^{t}F^{''}_{n}(B_{s})d{\langle
B \rangle _{s}}\\
&\longrightarrow F(B_{t})-F(0)+\frac12\int_{\mathbb R} f(x){\mathscr
L}(dx,t)
\end{align*}
in ${\mathbb L}^1$, as $n\to \infty$, and the theorem follows if $f$
is bounded uniformly in ${\mathbb R}$.

{\bf Step II.} Let $f$ be of bounded $p$-variation with $1\leq p<2$
and let $F$ be an absolutely continuous function with derivative
$F'=f$. Then $f$ is locally bounded on ${\mathbb R}$. Consider the
sets
$$
\Omega_k=\left\{\sup_{0\leq s\leq t}|B_s|<k\right\},\qquad
k=1,2,\ldots.
$$
Let $f^{[k]}\in {\mathcal W}_p$ be an uniformly bounded function on
${\mathbb R}$ such that $f^{[k]}=f$ on $[-k,k]$ and vanishes
outside, and let $F^{[k]}$ be an absolutely continuous function with
derivative $\frac{d}{dx}F^{[k]}=f^{[k]}$. Then the {\bf Step I}
implies that the formula
$$
F^{[k]}(B_t)=F^{[k]}(0)+\int_0^t
f^{[k]}(B_s)dB_s-\frac{1}{2}\int_{\mathbb R}f^{[k]}(x){\mathscr
L}(dx,t)
$$
holds q.s. on the set $\Omega_k$ for all $k=1,2\ldots$. Letting $k$
tend to infinity we deduce the desired It\^o
formula~\eqref{th3.1-eq1} if $f$ is of bounded $p$-variation with
$1\leq p<2$.
\end{proof}

\section{The quadratic covariation}\label{sec4}
In this section, we consider the quadratic covariation of $f(B)$ and
$B$, where $B$ is a $G$-Brownian motion and $f$ is a Borel function
on ${\mathbb R}$.
\begin{definition}
For a Borel function $f: \mathbb{R}\mapsto\mathbb{R}$ we define the
quadratic covariation $\langle f(B),B\rangle$ of $f(B)$ and $B$ as
follows
\begin{align*}
\langle f(B),B\rangle_{t}:=\lim_{n\rightarrow\infty}
\sum_{k=0}^{n-1}\{f(B_{t_{k+1}})-f(B_{t_{k}})\}(B_{t_{k+1}}-B_{t_{k}})
\end{align*}
in ${\mathbb L}^{1}$, where $\{t_{k}\}$ is a partitions of $[0,T]$
such that $\max\limits_{k}\{t_{k}-t_{k-1}\}\to 0$ as $n\to \infty$.
\end{definition}

\begin{lemma}\label{lem4.1}
For $f\in C^{1}(\mathbb{R})$, we
have
\begin{align*}
\langle f(B),B\rangle_{t}=\int_{0}^{t}f'(B_{u})d{\langle B \rangle
_{u}}.
\end{align*}
\end{lemma}
\begin{proof}
By the Taylor's expansion, we have
\begin{align*}
f(B_{t_{k+1}})-f(B_{t_{k}})=f'(B_{t_{k}})(B_{t_{k+1}}-B_{t_{k}})
+o(B_{t_{k+1}}-B_{t_{k}})
\end{align*}
for all $k=1,2,\ldots$. We claim that the convergence
\begin{align}\label{sec4-eq4.1}
\sum_{k=0}^{n-1}\left[o(B_{t_{k+1}}-B_{t_{k}})\right]
(B_{t_{k+1}}-B_{t_{k}})\longrightarrow 0
\end{align}
and
\begin{align}\label{sec4-eq4.2}
\sum_{k=0}^{n-1}f'(B_{t_{k}})(B_{t_{k+1}}-B_{t_{k}})^2
\longrightarrow \int_{0}^{t}f'(B_{u})d{\langle B \rangle _{u}}
\end{align}
hold in ${\mathbb L}^{1}$, as $n\to \infty$. The
convergence~\eqref{sec4-eq4.1} follows from the H\"older continuity
of $G$-Brownian motion. In order to prove~\eqref{sec4-eq4.2}, noting
that $f'(B_t)\in {\mathbb L}^1_G$, we get
\begin{equation}\label{sec4-eq4.3}
\hat{\mathbb{E}}\left[f'(B_s)f'(B_{s'})\int_{s}^{t}(B_{u}-B_{s})dB_{u}
\int_{s'}^{t'}(B_{u}-B_{s})dB_{u}\right]=0
\end{equation}
for $0<s'<t'<s<t$. In fact by
$$
f'(B_{s'})\int_{s'}^{t'}(B_{u}-B_{s})dB_{u}\in {\mathbb L}_G^1
$$
for $0<s'<t'<s<t$, we have
\begin{align*}
\hat{\mathbb{E}}&\left[f'(B_s)f'(B_{s'})\int_{s}^{t}(B_{u}-B_{s})dB_{u}
\int_{s'}^{t'}(B_{u}-B_{s})dB_{u}\right]\\
&
=\hat{\mathbb{E}}\left[\left(f'(B_{s'})\int_{s'}^{t'}(B_{u}-B_{s})dB_{u}
\right)^{+}\hat{\mathbb{E}}\left[\int_{s}^{t}(B_{u}-B_{s})dB_{u}
|\mathcal {F}_{s}\right]\right.\\
&\qquad+\left.\left(f'(B_{s'})\int_{s'}^{t'}(B_{u}-B_{s})dB_{u}
\right)^{-}\hat{\mathbb{E}}\left[-\int_{s}^{t}(B_{u}-B_{s})dB_{u}
|\mathcal {F}_{s} \right]\right]=0
\end{align*}
for $0<s'<t'<s<t$. On the other hand, we have
\begin{align*}
(B_{t_{k+1}}&-B_{t_{k}})^{2}=\langle B \rangle _{t_{k+1}}-\langle B
\rangle _{t_{k}}+2\int_{t_{k}}^{t_{k+1}}(B_{u}-B_{t_{k}})dB_{u}
\end{align*}
by using the identity
$B_{t}^{2}=\langle B \rangle _{t}+2\int_{0}^{t}B_{s}dB_{s}$.
Combining this with~\eqref{sec4-eq4.3}, we get
\begin{align*}
\hat{\mathbb{E}}&\left[\left|\sum_{k=0}^{n-1}f'(B_{t_{k}})
(B_{t_{k+1}}-B_{t_{k}})^{2}-\sum_{k=0}^{n-1}f'(B_{t_{k}})(\langle
B \rangle _{t_{k+1}}-\langle B \rangle _{t_{k}} )\right|^{2}\right]\\
&=\hat{\mathbb{E}}\left[\left|\sum_{k=0}^{n-1}2f'(B_{t_{k}})
\int_{t_{k}}^{t_{k+1}}(B_{u}-B_{t_{k}})dB_{u}\right|^{2}\right]\\
&\leq 4\sum_{k=0}^{n-1}\hat{\mathbb{E}}\left(f'(B_{t_{k}})
\int_{t_{k}}^{t_{k+1}}(B_{u}-B_{t_{k}})dB_{u}\right)^{2}\\
& =4\sum_{k=0}^{n-1}\hat{\mathbb{E}}\left[\left|f'(B_{t_{k}})
\right|^{2}
\right]\hat{\mathbb{E}}\left[\left(\int_{t_{k}}^{t_{k+1}}
(B_{u}-B_{t_{k}}
)dB_{u}\right)^{2}\right]\\
&\leq 4M^2\sum_{k=0}^{n-1}
\hat{\mathbb{E}}\left[\left(\int_{t_{k}}^{t_{k+1}} (B_{u}-B_{t_{k}})
dB_{u}\right)^{2}\right]\\
&\leq
4(M\bar{\sigma})^2\sum_{k=0}^{n-1}\int_{t_{k}}^{t_{k+1}}(u-t_{k})du=
(M\bar{\sigma})^2\sum_{k=0}^{n-1}(t_{k+1}-t_{k})^{2}\longrightarrow
0,
\end{align*}
as $n \to\infty $, where $M=\sup\limits_{0\leq s\leq t}\hat{\mathbb{E}}\left[\left|f'(B_s)\right|^{2}
\right]$, which deduces the convergence~\eqref{sec4-eq4.2}, and the lemma follows.
\end{proof}

As a direct consequence of the above lemma, we can rewrite the
It\^o's formula as follows.
\begin{corollary}\label{cor4.1}
If $F\in C^{2}({\mathbb R})$ and $F'=f$, we have
\begin{align*}
F(B_{t})=F(B_{0})+\int_{0}^{t}f(B_{u})dB_{u}+\frac{1}{2}\langle
f(B),B\rangle_{t},
\end{align*}
and the Bouleau-Yor identity
$$
\langle f(B),B\rangle_t=-\int_{\mathbb R}f(x){\mathscr
L}(dx,t),\qquad t\geq 0,
$$
holds.
\end{corollary}

\begin{theorem}\label{th4.1}
Let the function $x\mapsto f(x) $ be of bounded $p$-variation with
$1\leq p <2$ and let
\begin{equation}\label{sec4-eq4.4-00=00}
\begin{split}
\Bigl|\hat{\mathbb{E}}[f_z(B_t)f_z(B_{t'})&
(B_t-B_s)(B_{t'}-B_{s'})]\Bigr|\\
&\leq C\frac{(t-s)|t'-s'|}{t'(t-t')}\left|
\hat{\mathbb{E}}\left[f_z(B_t)f_z(B_{t'})B_{t'}(B_t-B_{t'})
\right]\right|
\end{split}
\end{equation}
for all $t>s\geq t',s'>0$ and $z\in {\mathbb R}$, where
$f_z(\cdot)=f(\cdot-z)-f(\cdot)$ for all $z\in {\mathbb R}$. Then
the quadratic covariation $\langle f(B),B\rangle$ exits in ${\mathbb
L}^1$, and the Bouleau-Yor identity
\begin{equation}\label{sec4-eq4.4}
\langle f(B),B\rangle_{t} =-\int_{\mathbb R}f(x){\mathscr L}(dx,t)
\end{equation}
holds q.s., for all $t\in[0,T]$.
\end{theorem}

It is important to note that the condition~\eqref{sec4-eq4.4-00=00}
is clear and it is an identity with $C=1$ if $B=W$ is a classical
standard Brownian motion. In order to see that the
condition~\eqref{sec4-eq4.4-00=00} holds for a classical standard
Brownian motion $W$, by approximating we may assume that $f$ is an
infinitely differentiable function with compact support and $t>s\geq
t'>s'>0$. It follows from the formula of integration by parts
(between divergence integral and Malliavin derivative operator $D$)
$$
E\left[F\int_0^Tu_sdW_s\right]=E\left[\int_0^TD_sFu_sds\right]
$$
with $F$ being a smooth random variable that
\begin{align*}
E&\left[f_z(W_{
t})f_z(W_{t'})(W_t-W_s)(W_{t'}-W_{s'}) \right]\\
\notag &=E\left[f_z(W_{
t})f_z(W_{t'})(W_t-W_s)\int_{s'}^{t'}dW_l\right]\\
&=E\int_0^T1_{[s',t']}(u)D_u\left[ f_z(W_{
t})f_z(W_{t'})(W_t-W_s)\right]du\\
&=(t'-s')E\left[ f_z'(W_{ t})f_z(W_{t'})(W_t-W_s)\right]\\
&\qquad+(t'-s')E\left[ f_z(W_{
t})f_z'(W_{t'})(W_t-W_s)\right]\\
&=(t'-s')(t-s)E\left[ f_z''(W_{
t})f_z(W_{t'})\right]+(t'-s')(t-s)E\left[ f'_z(W_{
t})f'_z(W_{t'})\right].
\end{align*}
For Malliavin calculus, see Nualart~\cite{Nualart2}. Let now
$$
\varphi_{t,s}(x,y)=\frac1{2\pi\rho_{t,s}}\exp\left(
-\frac1{\rho_{t,s}^2}\left(sx^2-2\mu_{t,s} xy+ty^2\right)\right)
$$
be the density function of $(W_t,W_s)$ with
$\mu_{t,s}=E(W_tW_s)=s\wedge t=s$ and
$\rho_{t,s}^2=st-\mu^2_{t,s}=s(t-s)$. We then have, by integration
by parts,
\begin{align*}
E\left[ f_z''(W_{ t})f_z(W_{t'})\right]&=\int_{\mathbb{R}^2}
f_z(x)f_z(y)\frac{\partial^2}{\partial x^2}\varphi_{t,t'}(x,y)dxdy\\
&=\int_{\mathbb{R}^2} f_z(x)f_z(y)\left\{\frac1{\rho^4_{t,t'}}
({t'}x-t'
y)^2-\frac{t'}{\rho^2_{t,t'}}\right\}\varphi_{t,t'}(x,y)dxdy
\end{align*}
\begin{align*}
E\left[ f'_z(W_{ t})f'_z(W_{t'})\right]&=\int_{\mathbb{R}^2}
f_z(x)f_z(y)\left\{\frac{1}{\rho^4_{t,t'}}(ty-t'x)(t'x-
t'y)+\frac{t'}{\rho^2_{t,t'}}\right\}\varphi_{t,t'}(x,y)dxdy.
\end{align*}
Combining these with the fact
$$
(ty-t'x)(t'x-t'y)=\rho^2_{t,t'}y(x-y)-(t')^2(x-y)^2,
$$
we get
\begin{align*}
E\bigl[ f_z''(W_{ t})f_z(W_{t'})&+f'_z(W_{ t})f'_z(W_{t'})\bigr]\\
&=\frac1{\rho^2_{t,t'}}\int_{\mathbb{R}^2}
f_z(x)f_z(y)(x-y)y\varphi_{t,t'}(x,y)dxdy\\
&=\frac1{\rho^2_{t,t'}}E\left[f_z(W_t)f_z(W_{t'})(W_t-W_{t'})W_{t'}
\right],
\end{align*}
which gives
\begin{align*} E&\left[f_z(W_{
t})f_z(W_{t'})(W_t-W_s)(W_{t'}-W_{s'})\right]\\
&=(t'-s')(t-s)
\frac1{\rho^2_{t,t'}}E\left[f_z(W_t)f_z(W_{t'})(W_t-W_{t'})W_{t'}
\right].
\end{align*}

Thus, the assumption~\eqref{sec4-eq4.4-00=00} is a natural condition
under the sublinear expectation $\hat{{\mathbb E}}$.

\begin{lemma}[Lemma 4.3 in
Dudley-Norvai\v{s}a~\cite{Dudley}]\label{lemma3.4} For $1\leq
p<\infty$, $f$ is of bounded $p$-variation on $[a,b]$ if and only if
$$
f=g\circ h
$$
for a bounded nondecreasing nonnegative function $h$ on $[a,b]$ and
a function on $[h(a),h(b)]$ satisfying a H\"older condition with
exponent $1/p$.
\end{lemma}

\begin{proof}[Proof of Theorem~\ref{th4.1}]
If $f\in C^{1}({\mathbb R})$, we then have by Corollary~\ref{cor4.1}
\begin{align*}
\langle f(B),B\rangle_{t}=&\int_{0}^{t}f'(B_{u})d{\langle B \rangle
_{u}} =-\int_{\mathbb R}f(x){\mathscr L}(dx,t).
\end{align*}
Let now $f \notin C^{1}({\mathbb R})$. By the localization argument
similar to proof of Theorem~\ref{th3.1} we may assume that {\em $f$
is uniformly bounded in the next discussion}. For $n\geq 1$ ,we set
\begin{align*}
f_{n}(x):=\int_{\mathbb R}\theta_{n}(x-y)f(y)dy,
\end{align*}
where $\theta_{n},\;n\geq 1$ are the mollifiers defined by~\eqref{sec3-eq3.7}. Then $f_{n}\in C^{\infty}({\mathbb R})$ is bounded and
\begin{equation}\label{sec4-eq4.5}
\langle f_{n}(B),B\rangle_{t} =-\int_{\mathbb R}f_{n}(x){\mathscr
L}(dx,t),\quad n\geq 1
\end{equation}
for all $t\in[0,T]$. Consider now the double sequence
\begin{align*}
\alpha_{m,t}(f_n)=\sum_{k=0}^{m-1}\left\{f_{n}(B_{t_{k+1}})
-f_{n}(B_{t_k})\right\}\left(B_{t_{k+1}}- B_{k}\right), \qquad
m,n\geq 1.
\end{align*}
In order to prove~\eqref{sec4-eq4.4}, we claim that
\begin{align*}
\hat{\mathbb{E}}\left[\left|\alpha_{m,t}(f) +\int_{\mathbb
R}f(x){\mathscr L}(dx,t)\right|\right] \longrightarrow 0,
\end{align*}
as $m\to \infty$. We have
\begin{equation}\label{eq4.5000}
\begin{split}
\hat{\mathbb{E}}\left[\left|\alpha_{m,t}(f)
+\int_{\mathbb R}f(x){\mathscr L}(dx,t)\right|\right]&\leq \hat{\mathbb{E}}\left[\left|\alpha_{m,t}(f)
-\alpha_{m,t}(f_n)\right|\right]\\
&\hspace{-1cm}+\hat{\mathbb{E}} \left[\left|\alpha_{m,t}(f_n)
+\int_{\mathbb
R}f_{n}(x){\mathscr L}(dx,t)\right|\right]\\
&\hspace{-1cm}+\hat{\mathbb{E}}\left[\left|\int_{\mathbb R}
f(x){\mathscr L}(dx,t)-\int_{\mathbb R} f_{n}(x){\mathscr
L}(dx,t)\right|\right].
\end{split}
\end{equation}
Let us estimate the three terms above. Denote
$\Delta_jB=B_{t_j}-B_{t_{j-1}}$ for $j=1,2,\ldots,m$ and
\begin{align*}
g_n(y)&=\int_{\mathbb R}\left[f(x-y)-f(y)\right]\theta_n(x)dx\\
&=\int_0^2\left[f(y-\frac{x}n)-f(y)\right]\theta(x)dx,\quad
n=1,2,\ldots.
\end{align*}
It follows that
\begin{align*}
\hat{\mathbb{E}}\left[|\alpha_{m,t}(f)-\alpha_{m,t}(f_n)| \right]^2
&=\hat{\mathbb{E}}\left( \sum_{j=1}^m(g_n(B_{t_j})-g_n(B_{t_{j-1}}))
\Delta_jB \right)^2\\
&\leq  \sum_{i<j}\hat{\mathbb{E}}(g_n(B_{t_j})-g_n(B_{t_{j-1}}))
(g_n(B_{t_i})-g_n(B_{t_{i-1}}))\Delta_jB\Delta_iB\\
&\qquad+\sum_{j=1}^m\hat{\mathbb{E}} (g_n(B_{t_j})
-g_n(B_{t_{j-1}}))^2(\Delta_jB)^2\\
&\leq \sum_{i<j}\hat{\mathbb{E}}(g_n(B_{t_j})-g_n(B_{t_{j-1}}))
(g_n(B_{t_i})-g_n(B_{t_{i-1}}))\Delta_jB\Delta_iB\\
&\qquad+\sum_{j=1}^m\hat{\mathbb{E}} \left[\left(g_n^2(B_{t_j})
+g_n^2(B_{t_{j-1}})\right) (\Delta_jB)^2\right]\\
&\equiv {\rm I+II}
\end{align*}
for all $t\geq 0$. Denote
\begin{equation}\label{sec4-eq4.4-01000}
\Delta_{n}(s,r,x,y):=\hat{\mathbb{E}}\left[(f(B_s-\frac{x}n)
-f(B_s))^2 (f(B_r-\frac{y}n)-f(B_r))^2\right]
\end{equation}
for any $s,r>0,x,y\in {\mathbb R}$ and $n=1,2,\ldots$. Notice that
\begin{align*}
\left|f(B_s-\frac{x}n)-f(B_s)\right|\leq 2\sup_x|f(x)|
\end{align*}
for all $s\geq 0$ and $n=1,2,\ldots$. We have, by the
condition~\eqref{sec4-eq4.4-00=00}
\begin{align*} |{\rm I}|&\leq \sum_{1\leq
i<j\leq m}\left|\hat{\mathbb{E}} \left[g_n(B_{t_j})g_n(B_{t_i})
\Delta_jB\Delta_iB\right] +\hat{\mathbb{E}}
\left[g_n(B_{t_j})g_n(B_{t_{i-1}}) \Delta_jB\Delta_iB\right]\right|\\
&\leq C\sum_{1\leq i<j\leq m}\frac{(t_j-t_{j-1})(t_i-t_{i-1})
}{\sqrt{t_i(t_j-t_i)}}\int_0^2\int_0^2\Delta_{n}(t_j,t_i,x,y)
\theta(x)\theta(y)dxdy\\
&\qquad+C\sum_{1\leq i<j\leq m}\frac{(t_j-t_{j-1})(t_i-t_{i-1})
}{\sqrt{t_{i-1}(t_j-t_{i-1})}}\int_0^2\int_0^2
\Delta_{n}(t_j,t_{i-1},x,y)\theta(x)\theta(y)dxdy\\
&\longrightarrow
C\int_0^tds\int_0^s\frac{dr}{\sqrt{r(s-r)}}\int_0^2\int_0^2
\Delta_{n}(s,r,x,y)\theta(x)\theta(y)dxdy
\end{align*}
and
\begin{align*}
|{\rm II}|&\leq C\overline{\sigma}^2 \sum_{j=1}^m
(t_j-t_{j-1})\int_0^2\Delta_{n}(t_j,t_j,x,x)\theta(x)dx\\
&\qquad+C\overline{\sigma}^2 \sum_{j=1}^m (t_j-t_{j-1})\int_0^2
\hat{{\mathbb E}}\left[(f(B_{t_{j-1}}+\frac{x}n)-f(B_{t_{j-1}}))^2
\right]\theta(x)dx\\
& \longrightarrow C\overline{\sigma}^2\int_0^tds\int_0^2
\Delta_{n}(s,s,x,x)\theta(x)dx\\
&\qquad\qquad+C\overline{\sigma}^2 \int_0^tds\int_0^2\hat{{\mathbb
E}}\left[(f(B_{s}+\frac{x}n)-f(B_{s}))^2\right]\theta(x)dx,
\end{align*}
as $m\to \infty$ for all $n\geq 1$. Thanks to Lemma~\ref{lemma3.4},
we see that there exists a bounded nondecreasing nonnegative
function $h$ such that
\begin{align*}
\left|f(B_s-\frac{x}n)-f(B_s)\right| &\leq
C_p\left(h(B_s)-h(B_s-\frac{x}n)\right)^{1/p}\downarrow 0
\end{align*}
q.s., as $n\to \infty$, which deduces
\begin{align*}
\lim_{n\to \infty}\lim_{m\to
\infty}\hat{\mathbb{E}}\left[\left|\alpha_{m,t}(f)
-\alpha_{m,t}(f_n)\right|\right]=0
\end{align*}
for all $t\geq 0$ by Lemma~\ref{lemma2.1000} and Lebesgue's
dominated convergence theorem. Noting that
$$
\langle f_n(B),B\rangle_t=-\int_{\mathbb R}f_n(x){\mathscr L}(dx,t)
$$
leads to
\begin{align*}
\hat{\mathbb{E}}\left[\left|\alpha_{m,t}(f_n) +\int_{\mathbb
R}f_{n}(x){\mathscr L}(dx,t)\right|\right] =
\hat{\mathbb{E}}\left[|\alpha_{m,t}(f_n)-\langle
f_n(B),B\rangle_t|\right]\longrightarrow 0,
\end{align*}
as $m\to \infty$, for all $n\geq 1$. It follows
from~\eqref{eq4.5000} and Lemma~\ref{lem3.3} that
\begin{align*}
\lim_{m\to \infty}\hat{\mathbb{E}}\left[\left|\alpha_{m,t}(f)
+\int_{\mathbb R}f(x){\mathscr L}(dx,t)\right|\right]&\leq
\lim_{n\to \infty}\lim_{m\to
\infty}\hat{\mathbb{E}}\left[\left|\alpha_{m,t}(f)
-\alpha_{m,t}(f_n)\right|\right]\\
&\hspace{-4cm}+\lim_{n\to \infty}\lim_{m\to \infty}\hat{\mathbb{E}}
\left[\left|\alpha_{m,t}(f_n) +\int_{\mathbb R}f_{n}(x){\mathscr
L}(dx,t)\right|\right]\\
&\hspace{-4cm}+\lim_{n\to
\infty}\hat{\mathbb{E}}\left[\left|\int_{\mathbb R} f(x){\mathscr
L}(dx,t)-\int_{\mathbb R} f_{n}(x){\mathscr L}(dx,t)\right|\right]
=0,
\end{align*}
and the theorem follows.
\end{proof}

According to Theorem~\ref{th3.1} and Theorem~\ref{th4.1}, we get an
extension of It\^{o}'s formula (F\"ollmer-Protter-Shiryayev's
formula).
\begin{corollary}
Let the function $f: \mathbb{R}\to \mathbb{R}$ be of bounded
$p$-variation with $1\leq p <2$ and let F be an absolutely
continuous function with derivative $F'=f$. If the
condition~\eqref{sec4-eq4.4-00=00} holds, we then have
\begin{align*}
F(B_{t})=F(0)+\int_{0}^{t}f(B_{s})dB_{s}+\frac{1}{2}\langle
f(B),B\rangle_{t}
\end{align*}
for all $t\in[0,T]$.
\end{corollary}
\section{Two parameter integrals of local time}\label{sec5}

In this section we turn to consider two parameter integrals
\begin{equation}\label{sec5-eq5.00}
\int_{0}^{t}\int_{\mathbb R}g(x,s){\mathscr L}(dx,ds).
\end{equation}

We first give an extension to the Lemma~\ref{lem3.1}.
\begin{lemma}[Occupation times formula]\label{lem5.1}
For all $t\geq 0$ and every bounded continuous Borel function
$\Phi$, we have
\begin{equation}\label{sec5-eq5.1}
\int_{0}^{t}\Phi(B_{s},s)d\langle
B\rangle_{s}=\int_{\mathbb R}dx\int_{0}^{t}\Phi(x,s){\mathscr
L}(x,ds) \qquad q.s.
\end{equation}
\end{lemma}
\begin{proof}
Let $\Phi(x,s)=1_{[a,b)}(x)1_{[t_1,t_2)}(s)$ with $a,b\in {\mathbb R}$ and $t_1,t_2\in [0,t]$. We have
\begin{align*}
\int_{0}^{t}\Phi(B_{s},s)d\langle
B\rangle_{s}=\int_{t_1}^{t_2}1_{[a,b)}(B_{s})d\langle
B\rangle_{s}
\end{align*}
and
\begin{align*}
\int_{\mathbb R}dx\int_{0}^{t}\Phi(x,s){\mathscr
L}(x,ds)&=\int_a^bdx\int_{t_1}^{t_2}{\mathscr L}(x,ds)\\
& =\int_a^bdx[{\mathscr
L}(x,t_2)-{\mathscr L}(x,t_1)]\\&
=\int_a^b{\mathscr
L}(x,t_2)dx-\int_a^b{\mathscr L}(x,t_1)dx\\
&=\int_{t_1}^{t_2}1_{[a,b)}(B_{s})d\langle
B\rangle_{s}\qquad q.s.
\end{align*}
which imply that~\eqref{sec5-eq5.1} holds. Consequently,~\eqref{sec5-eq5.1} deduces for the simple functions
$$
\Phi_{\Delta}(x,s)=\sum_{i,j}a_{ij} 1_{[x_{i},x_{i+1})}(x)1_{[t_{j},t_{j+1})}(s),
$$
where $\{a=x_{0}<x_{1}<\cdots<x_{n}=b, 0=t_{0}<x_{1}<\cdots<t_{n}=T\}$ is an arbitrary partition of $[a,b]\times[0,T]$.

On the other hand, for every bounded continuous Borel function
$\Phi(x,s)$ there is a sequence of bounded simple functions
$\Phi_{n}(x,s);n\geq 1$ such that $\Phi_{n}\to \Phi$ uniformly, in
${\mathbb R}\times {\mathbb R}_{+}$, and hence
\begin{align*}
\int_{0}^{t}\Phi_{n}(B_{s},s)d\langle B\rangle_{s}\longrightarrow
\int_{0}^{t}\Phi(B_{s},s)d{\langle B\rangle_{s}}
\end{align*}
in ${\mathbb L}^1$. It follows that there exists a subsequence
$\Phi_{n_{k}},k\geq 1$ such that
\begin{equation}
\int_{0}^{t}\Phi_{n_k}(B_{s},s)d\langle B\rangle_{s}\rightarrow
\int_{0}^{t}\Phi(B_{s},s)d\langle B\rangle_{s},\qquad q.s.
\end{equation}
and
\begin{equation}
\int_{\mathbb R}dx\int_{0}^{t}\Phi_{n_k}(x,s){\mathscr
L}(x,ds)\rightarrow \int_{\mathbb R}dx\int_{0}^{t}\Phi(x,s){\mathscr
L}(x,ds) ,\qquad q.s.
\end{equation}
which deduce
\begin{align*}  \int_{0}^{t}\Phi(B_{s},s)d\langle
B\rangle_{s}=\int_{\mathbb R}dx\int_{0}^{t}\Phi(x,s){\mathscr
L}(x,ds) \qquad q.s.
\end{align*}
and the lemma follows.
\end{proof}

In order to define the two parameter integrals~\eqref{sec5-eq5.00}
we use the idea from Feng-Zhao~\cite{Feng}. Recall that a function
$(x,y)\mapsto F(x,y)$, defined on $[a, b]\times [c, d]$ is of
bounded $p$-variation in $x$ uniformly in $y$, if
\begin{align*}
\sup_{y\in [c, d]}\sum
_{i=1}^{m}|F(x_{i},y)-F(x_{i-1},y)|^{p}<\infty,
\end{align*}
where $\{a=x_{0}<x_{1}<\cdot\cdot\cdot<x_{m}=b\}$ is an arbitrary
partition of $[a, b]$, and furthermore, it is of bounded
$p,q$-variation in $(x, y)$, if
\begin{align*}
\sup_{[a, b]\times [c, d]}\sum _{j=1}^{n}(\sum_{i=1}^{m}|\Delta
F(x_{i},y_{j})|^{p})^{q}<\infty,
\end{align*}
where
\begin{align*}
\Delta F(x_{i},y_{j})=F(x_{i},y_{j})-F(x_{i-1},y_{j})
-F(x_{i},y_{j-1}) +F(x_{i-1},y_{j-1}),
\end{align*}
and $\{a=x_{0}<x_{1}<\cdot\cdot\cdot<x_{m}=b;
c=y_{0}<y_{1}<\cdot\cdot\cdot<y_{n}=d\}$ is arbitrary partition of
$[a, b]\times [c, d]$.

Let now $G: \mathbb{R}\times \mathbb{R}\rightarrow \mathbb{R}$ be a
continuous function of bounded $q_{1}$-variation in $x$ uniformly in
$y$, and be of bounded $q_{2}$-variation in $y$ uniformly in $x$;
the continuous function $F: \mathbb{R}\times \mathbb{R}\rightarrow
\mathbb{R}$ be of bounded $p,q$-variation in $(x,y)$, where
$p,q,q_{1},q_{2}\geq 1$. Then the Young integral (see Theorem 3.1 in
Feng-Zhao~\cite{Feng})
\begin{align*}
\int_{a}^{b}\int_{c}^{d}G(x,s)F(dx,dy):=\lim_{\Delta_{m,n\rightarrow
0}}\sum_{j=1}^{n}\sum_{i=1}^{m}G(x_{i-1},y_{j-1})\Delta
F(x_{i},y_{j})
\end{align*}
is well defined, where
$\Delta_{m,n}=\max\limits_{i,j}\{|(x_{i},y_{j})-(x_{i-1},y_{j-1})|\}$, if
there exits two monotone increasing functions
$\rho:\mathbb{R}\rightarrow\mathbb{R}_{+}$ and
$\sigma:\mathbb{R}\rightarrow\mathbb{R}_{+}$ such that
\begin{align*}
\sum_{n,m}\rho(\frac{1}{n^{1/q_{1}}})\sigma(\frac{1}{n^{1/q_{2}}})
\frac{1}{n^{1/p}m^{1/pq}}<\infty.
\end{align*}
More works for two-parameter $p,q$-variation path integrals can be
fund in Feng-Zhao~\cite{Feng}. By taking $\rho(u)=u^{\alpha}$ and
$\sigma(u)=u^{1-\alpha}$ with $\alpha\in (0,1)$, one can prove the
following.
\begin{proposition}\label{prop5.1}
Let the continuous function
$F:\mathbb{R}\times[0,t]\rightarrow\mathbb{R}$ be of bounded
$p,q$-variation in $(x,t)$, and of bounded $\gamma$-variation in $x$
uniformly in $t$. Assume that the conditions
\begin{equation}\label{sec5-eq5.4}
1\leq\gamma<2,\qquad pq<q+\frac{1}{2},\qquad p,q\geq 1
\end{equation}
holds, then the Young integral (Bochner integral) of two parameters
\begin{align*}
\int_{0}^{t}\int_{\mathbb R}{\mathscr L}(x,s)F(dx,ds)
\end{align*}
exists q.s., and
\begin{align*}
\int_{0}^{t}\int_{\mathbb R}F(x,s){\mathscr
L}(dx,ds)=\int_{0}^{t}\int_{\mathbb R}{\mathscr
L}(x,s)F(dx,ds)-\int_{\mathbb R}{\mathscr L}(x,t)F(dx,t)
\end{align*}
\end{proposition}

Consider the smooth approximation of the function $F$
\begin{align}\label{sec5-eq5.5}
F_n(x,s):=\int_{0}^{2}\int_{0}^{2}\theta(r)\theta(z) F(x-\frac{r}{n},
s-\frac{z}{n})drdz, \qquad n\geq1,
\end{align}
where $\theta$ are the mollifiers defined in~\eqref{sec3-eq3.7}.
Then under the conditions of Proposition~\ref{prop5.1}, we have
\begin{align*}
\int_{0}^{t}\int_{\mathbb R}F_n(x,s){\mathscr L}(dx,ds)\rightarrow
\int_{0}^{t}\int_{\mathbb R}F(x,s){\mathscr L}(dx,ds)\qquad q.s.
\end{align*}
as $n\rightarrow\infty$.

\begin{corollary}\label{cor5.1-}
Under the conditions of Proposition~\ref{prop5.1}, we have
\begin{align*}
\int_{0}^{t}\int_{\mathbb R}F(x,s){\mathscr L}(dx,ds)=-\int_{\mathbb
R}dx\int_{0}^{t}\frac{\partial }{\partial x}F(x,s){\mathscr L}(x,ds)
\end{align*}
provided $F\in C^{1,1}(\mathbb{R}\times[0,T])$.
\end{corollary}

\begin{theorem}\label{th5.1}
$F\in C^{1,1}(\mathbb{R}\times\mathbb{R}_{+})$. Suppose that the
function $(x,t)\mapsto\frac{\partial }{\partial x}F(x,t)$ is of
bounded $p, q$-variation in $(x,t)$, and of bounded
$\gamma$-variation in $x$ uniformly in $t$, where these parameters
$p, q,\gamma$ satisfy the condition~\eqref{sec5-eq5.4}. Assume
that
$$
(x,t)\mapsto F(x,t),\qquad(x,t)\mapsto\frac{\partial }{\partial
x}F(x,t),\qquad (x,t)\mapsto\frac{\partial }{\partial t}F(x,t)
$$
are uniformly continuous in ${\mathbb R}^2$. Then the following
It\^o formula holds:
\begin{align*}
F(B_{t},t)=F(0,0)&+\int_{0}^{t}\frac{\partial }{\partial
t}F(B_{s},s)ds\\
&+\int_{0}^{t}\frac{\partial }{\partial
x}F(B_{s},s)dB_{s}
-\frac{1}{2}\int_{0}^{t}\int_{\mathbb
R}\frac{\partial }{\partial x}F(x,s){\mathscr L}(dx,ds).
\end{align*}
\end{theorem}
\begin{proof}
Let $F_n$ be defined in~\eqref{sec5-eq5.5} for $n\geq1$. Then
$F_n\in C^{2,1}(\mathbb{R}\times\mathbb{R}_{+})$ and $\frac{\partial
}{\partial x}F_n,n\geq1$ are of bounded $p, q$-variation in $(x,t)$,
and of bounded $\gamma$-variation in $x$ uniformly in $t$, and we
have
\begin{align*}
F_n(B_{t},t)=F_n(0,0)&+\int_{0}^{t}\frac{\partial }{\partial
t}F_n(B_{s},s)ds\\
&+\int_{0}^{t}\frac{\partial }{\partial
x}F_n(B_{s},s)dB_{s}-\frac12\int_{0}^{t}\int_{\mathbb
R}\frac{\partial }{\partial x}F_n(x,s){\mathscr L}(dx,ds),
\end{align*}
for all $n\geq 1$. On the other hand, we have
\begin{align*}
\hat{{\mathbb E}}&\left|\int_{0}^{t}\frac{\partial }{\partial
t}F_n(B_{s},s)ds-\int_{0}^{t}\frac{\partial }{\partial
t}F_n(B_{s},s)ds\right|\\
&\leq \int_{0}^{t}\int_{0}^{2}\int_{0}^{2}\theta(r)\theta(z)
\hat{{\mathbb E}}\left|\frac{\partial }{\partial
t}F(B_s-\frac{r}{n}, s-\frac{z}{n})-\frac{\partial }{\partial
t}F(B_s,s)\right|drdz\longrightarrow 0,
\end{align*}
and
\begin{align*}
\hat{{\mathbb E}}&\left|\int_{0}^{t}\frac{\partial }{\partial
x}F_n(B_{s},s)dB_s-\int_{0}^{t}\frac{\partial }{\partial
x}F_n(B_{s},s)dB_s\right|^2\\
&\leq C\bar{\sigma}^2\hat{{\mathbb
E}}\int_{0}^{t}\left(\frac{\partial }{\partial
x}F_n(B_{s},s)dB_s-\int_{0}^{t}\frac{\partial }{\partial
x}F_n(B_{s},s)\right)^2ds\\
&\leq
C\bar{\sigma}^2\int_{0}^{t}\int_{0}^{2}\int_{0}^{2}\theta(r)\theta(z)
\hat{{\mathbb E}}\left|\frac{\partial }{\partial
x}F(B_s-\frac{r}{n}, s-\frac{z}{n})-\frac{\partial }{\partial
x}F(B_s,s)\right|^2drdz\longrightarrow 0,
\end{align*}
as $n\to \infty$. Moreover, the occupation times
formula~\eqref{sec5-eq5.1} and Corollary~\ref{cor5.1-} imply that
\begin{align*}
\int_{0}^{t}\frac{\partial ^{2}}{\partial
x^{2}}F_n(B_{s},s)d\langle B\rangle_{s}&=\int_{\mathbb
R}dx\int_{0}^{t}\frac{\partial ^{2}}{\partial
x^{2}}F_n(x,s){\mathscr L}(x,ds)\\
&=-\int_{0}^{t}\int_{\mathbb
R}\frac{\partial }{\partial x}F_n(x,s){\mathscr
L}(dx,ds)\\
&\rightarrow -\int_{0}^{t}\int_{\mathbb R}\frac{\partial
}{\partial x}F(x,s){\mathscr L}(dx,ds)
\end{align*}
q.s., as $n\to \infty$. Combining this with $F_n\longrightarrow F$
uniformly, as $n\to \infty$, we obtain the theorem.
\end{proof}
As the end of this paper, we consider the quadratic covariation
$\langle f(B,\cdot),B\rangle$ of $f(B,\cdot)$ and $B$ defined by
\begin{align*}
\langle f(B,\cdot),B \rangle:=\lim_{n\to \infty}\sum_{j=1}^n
\{f(B_{t_{j}},t_{j})-f(B_{t_{j-1}},t_{j-1})\} (B_{t_j}-B_{t_{j-1}})
\end{align*}
in ${\mathbb{L}}^{1}$, where $f: \mathbb{R}\times [0,\infty)\to
\mathbb{R}$ is a Borel function and $\{t_j\}$ is a partitions of
$[0,T]$ such that $|\Delta_n|:=\max\{t_j-t_{j-1}\}\to 0$, as $n\to
\infty$.
\begin{theorem}\label{th5.2}
Let $f\in C(\mathbb{R}\times\mathbb{R}_{+})$ be of bounded
$p,q$-variation in $(x,t)$, and of bounded $\gamma$-variation in $x$
uniformly in $t$, and let the parameters $p,q,\gamma$ satisfy the
condition~\eqref{sec5-eq5.4}. Assume that
\begin{equation}\label{sec4-eq4.4-00=00==0}
\begin{split}
\Bigl|\hat{\mathbb{E}}[f_z(B_t,t)f_z(B_{t'},t')&
(B_t-B_s)(B_{t'}-B_{s'})]\Bigr|\\
&\leq C\frac{(t-s)|t'-s'|}{t'(t-t')}\left|
\hat{\mathbb{E}}\left[f_z(B_t,t)f_z(B_{t'},t')B_{t'}(B_t-B_{t'})
\right]\right|
\end{split}
\end{equation}
for all $t>s\geq t',s'>0$ and $z\in {\mathbb R}$, where
$f_z(\cdot)=f(\cdot+z)-f(\cdot)$ for all $z\in {\mathbb R}$, if the
function $f$ is uniformly continuous in ${\mathbb R}^2$, then the
quadratic covariation $\langle f(B,\cdot),B\rangle$ exits in
${\mathbb L}^1$, and the Bouleau-Yor identity
\begin{equation}\label{sec5-eq5.6}
\langle f(B,\cdot),B \rangle_{t}=-\int_{0}^{t}\int_{\mathbb
R}f(x,s){\mathscr L}(dx,ds)
\end{equation}
holds q.s. for all $t\in[0,T]$.
\end{theorem}
\begin{proof}
Let $f\in C^{1,1}(\mathbb{R}\times \mathbb{R}_{+})$. Similar to the
proof of Lemma~\ref{lem4.1} we can get
\begin{align*}
\langle f(B,\cdot),B \rangle_{t}=\int_{0}^{t}\frac{\partial
}{\partial x}f(B_{s},s)d\langle B\rangle_{s}
\end{align*}
for all $t\in[0,T]$. It follows from the occupation times
formula~\eqref{sec5-eq5.1} that
\begin{align*}
\langle f(B,\cdot),B \rangle_{t}&=\int_{\mathbb
R}dx\int_{0}^{t}\frac{\partial }{\partial
x}f(x,s){\mathscr L}(x,ds)\\
&=-\int_{0}^{t}\int_{\mathbb R}f(x,s){\mathscr L}(dx,ds)
\end{align*}
for all $t\in[0,T]$.

Let now $f\not\in C^{2,1}(\mathbb{R}\times \mathbb{R}_{+})$. For
$n\geq1$ we define $f_n$ as follows
\begin{align*}
f_n(x,s):=\int_{0}^{2}\int_{0}^{2}\theta(r)\theta(z)
f(x-\frac{r}{n}, s-\frac{z}{n})drdz, \qquad n\geq1,
\end{align*}
where $\theta$ is the mollifier defined in~\eqref{sec3-eq3.7}. Then
$f_n\in C^{1,1}(\mathbb{R}\times[0,T])$ and we have
\begin{align*}
\langle f_n(B,\cdot),B \rangle_{t}=-\int_{0}^{t}\int_{\mathbb
R}f_n(x,s){\mathscr L}(dx,ds)
\end{align*}
for all $t\in[0,T]$. Thus, similar to the proof of Theorem~\ref{th4.1} we can obtain the identity~\eqref{sec5-eq5.6}.
\end{proof}

According to the above theorems, we get an analogue of
F\"{o}llmer-Protter-Shiryayev's formula (see~\cite{Follmer}).
\begin{corollary}
Let $F\in C^{1,1}(\mathbb{R}\times \mathbb{R}_{+})$and let the
conditions in Theorem~\ref{th5.1} and Theorem~\ref{th5.2} hold. Then
\begin{align*}
F(B_{t},t)=F(0,0)+\int_{0}^{t}\frac{\partial }{\partial
t}F(B_{s},s)ds+\int_{0}^{t}f(B_{s},s)dB_{s}+\frac{1}{2}\langle
f(B,\cdot),B \rangle_{t}
\end{align*}
for all $t\in [0,T]$.
\end{corollary}

%\appendix

\section{{\bf Appendix}: Quadratic variation}\label{sec6}

In this appendix, we give the other representation of quadratic
covariation $\langle B,B\rangle$. From the properties of quadratic
variation of $G$-Brownian motion we have
\begin{align*}
\frac1\varepsilon\hat{\mathbb{E}}\int_0^{\varepsilon}
\langle B\rangle_rdr\longrightarrow 0,
\end{align*}
as $\varepsilon$ tends to $0$.

\begin{lemma}\label{lem6-1}
For all $t\geq 0$, we have
\begin{align*}
\frac1\varepsilon\int_{0}^{t}(\langle B\rangle_{s+\varepsilon}-\langle B\rangle_{s})ds \longrightarrow \langle B\rangle_{t}
\end{align*}
in ${\mathbb L}^1$, as $\varepsilon$ tends to $0$.
\end{lemma}
\begin{proof}
Notice that
\begin{align*}
\frac1\varepsilon\int_{0}^{t}(\langle B\rangle_{s+\varepsilon}
-\langle B\rangle_{s})ds&=\frac1\varepsilon\int_t^{t+\varepsilon}
\langle B\rangle_rdr-\frac1\varepsilon\int_0^{\varepsilon}
\langle B\rangle_rdr\\
&=\frac1\varepsilon\int_0^{\varepsilon}
\langle B\rangle_{t+r}dr-\frac1\varepsilon\int_0^{\varepsilon}
\langle B\rangle_rdr.
\end{align*}
We get
\begin{align*}
\hat{\mathbb{E}}\left|\frac1\varepsilon\int_{0}^{t}(\langle B\rangle_{s+\varepsilon}-\langle B\rangle_{s})ds-\langle B\rangle_t\right|&\leq \hat{\mathbb{E}}\left|\frac1\varepsilon\int_{0}^\varepsilon(
\langle B\rangle_{t+r}-\langle B\rangle_t)dr\right|
+\hat{\mathbb{E}}\frac1\varepsilon\int_0^{\varepsilon}
\langle B\rangle_rdr\\
&\leq \hat{\mathbb{E}}\left|\frac1\varepsilon\int_{0}^\varepsilon
\langle B^t\rangle_rdr\right|
+\hat{\mathbb{E}}\frac1\varepsilon\int_0^{\varepsilon}
\langle B\rangle_rdr\\
&\longrightarrow 0
\end{align*}
for all $t\geq 0$, as $\varepsilon$ tends to $0$.
\end{proof}
From the proof and the definition of the integral with respect to
$\langle B\rangle$ we get
\begin{equation}\label{sec6-eq6.1}
\frac1\varepsilon\int_{0}^{t}g(B_s)(\langle B\rangle_{
s+\varepsilon}-\langle B\rangle_{s})ds \longrightarrow
\int_{0}^{t}g(B_s)d\langle B\rangle_s
\end{equation}
in ${\mathbb L}^1$, as $\varepsilon$ tends to $0$, for all $g\in
C({\mathbb R})$.

\begin{proposition}\label{lem3-1.1}
Let $g\in C({\mathbb R})$. Then, for all $t\geq 0$, we have
\begin{align*}
\frac1\varepsilon\int_{0}^{t}g(B_s)(B_{s+\varepsilon}-B_{s})^2ds \longrightarrow \int_{0}^{t}g(B_s)d\langle B\rangle_s
\end{align*}
in ${\mathbb L}^1$, as $\varepsilon$ tends to $0$.
\end{proposition}
\begin{proof}
It is enough to show that the following convergence holds:
\begin{align*}
\frac1\varepsilon\int_{0}^{t}(B_{s+\varepsilon}-B_{s})^2ds
\longrightarrow \int_{0}^{t}d\langle B\rangle_s=\langle B\rangle_t
\end{align*}
in ${\mathbb L}^1$, as $\varepsilon$ tends to $0$. Clearly, we have
\begin{align*}
(B_{s+\varepsilon}-B_{s})^2=\langle B\rangle_{s+\varepsilon}-\langle B\rangle_{s}+2\int_s^{s+\varepsilon}(B_r-B_s)dB_r
\end{align*}
for all $s\geq 0$. Thus, it is enough to show
\begin{equation}\label{sec6-eq6.2}
\hat{\mathbb{E}}\left|\frac1{\varepsilon}\int_0^tds \int_s^{s+\varepsilon}(B_r-B_s)dB_r\right|\longrightarrow 0
\end{equation}
for all $t\geq 0$, as $\varepsilon$ tends to $0$, by~\eqref{sec6-eq6.1}. To end this we have, by~\eqref{sec4-eq4.3},
\begin{align*}
\frac1{\varepsilon^2}\hat{\mathbb{E}}&\left|\int_0^tds \int_s^{s+\varepsilon}(B_r-B_s)dB_r\right|^2\\ &=\frac1{\varepsilon^2}
\hat{\mathbb{E}}\left[\int_0^t \int_0^t\int_s^{s+\varepsilon}(B_u-B_s)dB_u \int_l^{l+\varepsilon}(B_v-B_l)dB_vdsdl
\right]\\
&\leq \frac1{\varepsilon^2}
\left[\int_0^tds \int_{s-\varepsilon}^s\hat{\mathbb{E}} \int_s^{s+\varepsilon}(B_u-B_s)dB_u \int_l^{l+\varepsilon}(B_v-B_l)dB_vdl
\right]\\
&\quad+\frac1{\varepsilon^2}
\left[\int_0^tds \int_s^{s+\varepsilon}\hat{\mathbb{E}} \int_s^{s+\varepsilon}(B_u-B_s)dB_u \int_l^{l+\varepsilon}(B_v-B_l)dB_vdl
\right].
\end{align*}
Notice that
\begin{align*}
\hat{\mathbb{E}}&\left(\int_s^{s+\varepsilon}(B_u-B_s)dB_u \int_l^{l+\varepsilon}(B_v-B_l)dB_v\right)\\
&\leq \hat{\mathbb{E}}\left((\int_s^{l+\varepsilon}(B_u-B_s)dB_u)^2 \right)+
\hat{\mathbb{E}}\left((\int_s^{l+\varepsilon}(B_u-B_s)dB_u)
(B_s-B_l)(B_{l+\varepsilon}-B_s)\right)\\
&\leq C \bar{\sigma}^2\int_s^{l+\varepsilon}(u-s)du=C \bar{\sigma}^2(l+\varepsilon-s)^2 +C\bar{\sigma}^2(l+\varepsilon-s)^{\frac32}(s-l)^{\frac12}
\end{align*}
for $l<s<l+\varepsilon<s+\varepsilon$, and
\begin{align*}
\hat{\mathbb{E}}&\left(\int_s^{s+\varepsilon}(B_u-B_s)dB_u \int_l^{l+\varepsilon}(B_v-B_l)dB_v\right)\\
&\leq C \bar{\sigma}^2\int_l^{s+\varepsilon}(u-s)du=C \bar{\sigma}^2(s+\varepsilon-l)^2 +C\bar{\sigma}^2(s+\varepsilon-l)^{\frac32}(l-s)^{\frac12}
\end{align*}
for $s<l<s+\varepsilon<l+\varepsilon$. We get
\begin{align*}
\frac1{\varepsilon^2}\hat{\mathbb{E}}&\left|\int_0^tds \int_s^{s+\varepsilon}(B_r-B_s)dB_r\right|^2\leq C\bar{\sigma}^2t\varepsilon\longrightarrow 0,
\end{align*}
as $\varepsilon$ tends to $0$, and the lemma follows.
\end{proof}

\begin{corollary}\label{cor3-1.1}
We have
\begin{align*}
\frac1\varepsilon\int_{0}^{t}(B_{s+\varepsilon}-B_{s})^2ds
\longrightarrow\langle B\rangle_{t}
\end{align*}
in ${\mathbb L}^1$, as $\varepsilon\to 0$.
\end{corollary}

\end{document}